\documentclass[a4paper,11pt,twoside]{article}
\usepackage{amssymb,amsmath,latexsym,geometry,fancyhdr,lineno,hyperref,titletoc}
\contentsmargin{0pt}

\dottedcontents{section}[2em]{\vspace{-1mm}\small}{2em}{0pt}
\dottedcontents{subsection}[5em]{\vspace{-1mm}\small}{3em}{5pt}

\hypersetup{colorlinks,linkcolor= blue,citecolor=blue}
\newtheorem{theorem}{Theorem}[section]

\newtheorem{lemma}[theorem]{Lemma}

\newtheorem{remark}[theorem]{Remark}

\newtheorem{definition}[theorem]{Definition}
\def\proof{\mbox {\textbf{Proof.}~~}}
\numberwithin{equation}{section}
\allowdisplaybreaks[0]

\begin{document}
\title{{\bf\Large On global minimizers for a mass constrained problem}}
\author{\\
{ \textbf{\normalsize Louis Jeanjean}}\\
{\it\small Laboratoire de Math\'{e}matiques (CNRS UMR 6623),}\\
{\it\small Universit\'{e} de Bourgogne Franche-Comt\'{e},}\\
{\it\small Besan\c{c}on 25030, France}\\
{\it\small e-mail: louis.jeanjean@univ-fcomte.fr}\\
\\
{ \textbf{\normalsize Sheng-Sen Lu}}\\
{\it\small LMAM and School of Mathematical Sciences,}\\
{\it\small Peking University,}\\
{\it\small Beijing 100871, PR China}\\
{\it\small e-mail: sslu@pku.edu.cn}}
\date{}
\maketitle
{\bf\normalsize Abstract.} {\small
In any dimension $N \geq 1$, for given mass $m > 0$ and for the $C^1$ energy functional
  \begin{equation*}
    I(u):=\frac{1}{2}\int_{\mathbb{R}^N}|\nabla u|^2dx-\int_{\mathbb{R}^N}F(u)dx,
  \end{equation*}
we revisit the classical problem of finding conditions on $F \in C^1(\mathbb{R},\mathbb{R})$ insuring that $I$ admits global minimizers on the mass constraint
  \begin{equation*}
    S_m:=\left\{u\in H^1(\mathbb{R}^N)~|~\|u\|^2_{L^2(\mathbb{R}^N)}=m\right\}.
  \end{equation*}
Under assumptions that we believe to be nearly optimal, in particular without assuming that $F$ is even, any such global minimizer, called energy ground state, proves to have constant sign and to be radially symmetric monotone with respect to some point in $\mathbb{R}^N$. Moreover, we show that any energy ground state is a least action solution of the associated action functional. This last result answers positively, under general assumptions, a long standing issue.}




\pagestyle{fancy}
\fancyhead{} 
\fancyfoot{} 
\renewcommand{\headrulewidth}{0pt}
\renewcommand{\footrulewidth}{0pt}
\fancyhead[CE]{ \textsc{Louis Jeanjean \& Sheng-Sen Lu}}
\fancyhead[CO]{ \textsc{On Global Minimizers for a Mass Constrained Problem}}
\fancyfoot[C]{\thepage}


\newpage
\section{Introduction}\label{sect:introduction}

Let $N \geq 1$ and $I: H^1(\mathbb{R}^N) \to \mathbb{R}$ be a $C^1$ energy functional defined by
  \begin{equation*}
    I(u):=\frac{1}{2}\int_{\mathbb{R}^N}|\nabla u|^2dx-\int_{\mathbb{R}^N}F(u)dx
  \end{equation*}
where $F(t) := \int^t_0 f(\tau) d\tau$ for some function $f \in C(\mathbb{R}, \mathbb{R})$.

In this paper we focus on the minimization problem
  \begin{equation}\tag{$\text{Inf}_m$}\label{eq:Inf_m}
    E_m := \inf_{u \in S_m} I(u),
  \end{equation}
where $m > 0$ is prescribed and
  \begin{equation*}
    S_m := \left\{u \in H^1(\mathbb{R}^N) ~ | ~ \|u\|^2_{L^2(\mathbb{R}^N)} = m\right\}.
  \end{equation*}
By a direct application of Lagrange multiplier's rule, if $u \in S_m$ solves \eqref{eq:Inf_m} then	there exists $\mu = \mu(u) \in \mathbb{R}$ such that
  \begin{equation}\label{eq:equation-libre}
    - \Delta u = f(u) - \mu u  \qquad \mbox{in} ~ H^1(\mathbb{R}^N).
  \end{equation}
A minimizer of \eqref{eq:Inf_m} is often called an energy ground state and $E_m$ the ground state energy.
		
The study of problem \eqref{eq:Inf_m}  naturally arises in the search of standing waves for nonlinear scalar field equations of the form
\begin{linenomath*}
    \begin{equation}\label{eq:equation-evolution}
     i \psi_t + \Delta \psi + f(\psi) =0, \qquad \psi : \mathbb{R}  \times \mathbb{R}^N \to \mathbb{C}.
    \end{equation}
  \end{linenomath*}
By standing waves, we mean solutions to \eqref{eq:equation-evolution} of the special form $\psi(t,x) = e^{i\mu t} u(x)$ with $\mu \in \mathbb{R}$ and $u \in H^1(\mathbb{R}^N)$. Clearly $\psi(t,x)$ satisfies \eqref{eq:equation-evolution} if $u(x)$ satisfies \eqref{eq:equation-libre} for the corresponding $\mu  \in \mathbb{R}$.

The study of such type of equations, which already saw major contributions forty years ago,  \cite{Be83-1, CL82, Lions84-1,Lions84-2, St82}, now lies at the root of several models linked with current physical applications (such as nonlinear optics, the theory of water waves, ...).   For these equations, finding solutions with a prescribed $L^2$-norm is particularly relevant since this quantity is preserved along the time evolution. In addition, if the solutions correspond to energy ground states, then, in most situations, it is possible to prove that the associated standing waves are orbitally stable. This likely explains why the study of problem \eqref{eq:Inf_m} is still the object of an intense activity. Among many others possible choices, we refer to \cite{CKS21, CS20, C03, DST21, HaSt04, HT19, JL19, LRN20, Sh14,  St19} and to the references therein.

Our first main result concerns the solvability of \eqref{eq:Inf_m}.  It encompasses several results previously proved for specific nonlinearities and, in particular, can be viewed as an extension of the one of \cite{Sh14} already obtained in a very general setting.  The following assumptions on  $f \in  C(\mathbb{R},\mathbb{R})$ will be required.
 \begin{itemize}
    \item[$(f1)$] \hypertarget{f1}{} $\lim_{t \to 0} f(t)/t = 0$.
		\item[$(f2)$] \hypertarget{f2}{} When $N \geq 3$,
                    \begin{equation*}
                      \limsup_{|t| \to \infty} \frac{|f(t)|}{|t|^{\frac{N + 2}{N - 2}}} < \infty;
                    \end{equation*}
                  when $N = 2$,
                    \begin{equation*}
                      \lim_{|t| \to \infty} \frac{f(t)}{e^{\alpha t^2}} = 0 \quad \forall \alpha > 0;
                    \end{equation*}
                  and also for any $N \geq 1$,
                    \begin{equation*}
                      \limsup_{|t| \to \infty} \frac{f(t)t}{|t|^{2 + \frac{4}{N}}} \leq 0.
                    \end{equation*}
    \item[$(f3)$] \hypertarget{f3}{} There exists $\zeta \neq 0$ such that $F(\zeta) > 0$.
  \end{itemize}

\begin{theorem}\label{thm:minimizers}
  Assume that $N \geq 1$, $f \in C(\mathbb{R}, \mathbb{R})$ satisfies $(\hyperlink{f1}{f1}) - (\hyperlink{f3}{f3})$. Then
    \begin{equation*}
      E_m := \inf_{u \in S_m} I(u) > - \infty
    \end{equation*}
  and the map $m \mapsto E_m$ is nonincreasing and continuous. Moreover,
	\begin{itemize}
      \item[$(i)$] there exists a number $m^* \in [0, \infty)$ such that
                     \begin{equation*}
                       E_m = 0 \quad \text{if} ~ 0 < m \leq m^*, \qquad E_m < 0 \quad \text{when} ~ m > m^*;
                     \end{equation*}
      \item[$(ii)$] when $m > m^*$, the global infimum $E_m$ is achieved and thus \eqref{eq:Inf_m} has an energy ground state $v \in S_m$ with $I(v) = E_m < 0$;
	  \item[$(iii)$] when $0 < m < m^*$, $E_m = 0$ is not  achieved;
	  \item[$(iv)$] $m^* = 0$ if in addition
                      \begin{equation*}\tag{A.1}\label{eq:key_1}
                        \lim_{t \to 0} \frac{F(t)}{|t|^{2 + \frac{4}{N}}} = + \infty,
                      \end{equation*}
                    and $m^* > 0$ if in addition
                      \begin{equation*}\tag{A.2}\label{eq:key_2}
                        \limsup_{t \to 0}\frac{F(t)}{|t|^{2 + \frac{4}{N}}} < + \infty.
                      \end{equation*}
    \end{itemize}
\end{theorem}
\begin{remark}\label{rmk:minimizers}
  \begin{itemize}
    \item[$(i)$] As it will be clear from the proof of Theorem \ref{thm:minimizers} $(ii)$, see also Remark \ref{rmk:compactness}, when $m > m^*$ we also show that any minimizing sequence for \eqref{eq:Inf_m} is, up to a subsequence and up to translations in $\mathbb{R}^N$, strongly convergent.
    \item[$(ii)$] When $0< m < m^*$, it is proved in Theorem \ref{thm:minimizers} $(iii)$ that the global infimum $E_m = 0$ is not achieved, but this does not mean that the constrained functional $I_{|S_m}$ may not admit critical points with positive energies, see the companion work \cite{JL21}.
    \item[$(iii)$] In the case $m^* > 0$, studying existence and nonexistence of global minimizers with respect to $E_{m^*} = 0$ seems to be a delicate issue. Since it exceeds our scope of the present paper, we shall not explore further general conditions on $f$ that ensure the existence or nonexistence but refer the interested reader to \cite[Lemma 2.3]{JL21} and \cite[Theorem 1.4]{Sh14} for some existence results.
    \item[$(iv)$] For convenience of statement, we introduce the notation
          \begin{equation*}
            m \succeq_f m^*
          \end{equation*}
        with the understanding that $m \geq m^*$ if $m^* > 0$ and $E_{m^*} = 0$ is achieved, and $m > m^*$ if otherwise. As one may observe, when $m \succeq_f m^*$ and for any minimizer $v \in S_m$ of \eqref{eq:Inf_m}, the associated Lagrange multiplier $\mu = \mu(v)$ is positive. Indeed, from the Pohozaev identity corresponding to \eqref{eq:equation-libre}, see \cite[Proposition 1]{Be83-1},
          \begin{equation*}
            P(v) := \frac{N - 2}{2N}\int_{\mathbb{R}^N} |\nabla v|^2 dx + \frac{1}{2} \mu \int_{\mathbb{R}^N} |v|^2 dx - \int_{\mathbb{R}^N} F(v) dx = 0
          \end{equation*}
        and the fact that $I(v) =  E_m \leq 0$, we have
          \begin{equation*}
            0 \geq I(v) = I(v) - P(v) = \frac{1}{N} \int_{\mathbb{R}^N} |\nabla v|^2 dx - \frac{1}{2} \mu m
          \end{equation*}
        and hence $\mu > 0$.
  \end{itemize}
\end{remark}

\begin{remark}\label{rmk:examples}
  Let us give some examples of nonlinearities satisfying $(\hyperlink{f1}{f1}) - (\hyperlink{f3}{f3})$.
    \begin{itemize}
      \item[$(i)$] $f(t) = |t|^{p-2}t + A |t|^{q-2}t$ with
                     \begin{equation*}
                       A \in \mathbb{R} \qquad \text{and} \qquad 2 < q < p < 2 + \frac{4}{N}.
                     \end{equation*}
                   In particular, \eqref{eq:key_1} and \eqref{eq:key_2} hold when $A \geq 0$ and when $A < 0$ respectively.
	  \item[$(ii)$] $f(t) = |t|^{p-2}t - |t|^{q-2}t$ with
                      \begin{equation*}
                        2 < p < q
                          \left\{
                            \begin{aligned}
                              & < \infty, & \qquad  & \text{if} ~ N = 1,2,\\
                              & \leq \frac{2N}{N-2}, & \qquad & \text{if} ~ N \geq 3.
                            \end{aligned}
                          \right.
                      \end{equation*}
                    In particular, \eqref{eq:key_1} and \eqref{eq:key_2} hold if $p < 2 + \frac{4}{N}$ and if $p \geq 2 + \frac{4}{N}$ respectively, and  when $N = 1,2,3$ we cover the so-called cubic-quintic nonlinearity
                      \begin{equation*}
                        f(t) = |t|^2t - |t|^4t
                      \end{equation*}
                    which attracts much attention due to its physical relevance, see for example \cite{CKS21,CS20,KOPV17,LRN20, St19}.
										
    \end{itemize}
 These examples  are only some special cases, our Theorem \ref{thm:minimizers} and the subsequent Theorems \ref{thm:properties} and \ref{thm:LeastAction} apply to more general nonlinearities, in particular to those which are not a sum of powers.
\end{remark}

The conditions $(\hyperlink{f1}{f1})$ and $(\hyperlink{f2}{f2})$ show that the energy functional $I$ is coercive and bounded from below on the sphere $S_m$. To prove the existence of minimizers in Theorem \ref{thm:minimizers}, instead of using directly the machinery of the Concentration Compactness Principle as presented in \cite{Lions84-1,Lions84-2}, we follow a version of it introduced in \cite{Ik14} which requires less technicalities. The vanishing scenario is ruled out by the fact that $E_m$ is negative for given $m > m^*$, and the dichotomy case is disproved by exploiting Lemma \ref{lemma:E_m} $(iv)$. To derive the sharp threshold mass $m^*$, we need the monotonicity and continuity of $E_m$ which is proved in Lemma \ref{lemma:E_m} $(v)$. As to the condition $(\hyperlink{f3}{f3})$, it plays its due role in Lemma \ref{lemma:E_m} $(ii)$ to show that $E_m$ is negative for large mass.

Our next result shows that any energy ground state has constant sign and enjoys symmetry and monotonicity properties. It unifies several more specific results scattered in the literature and in particular allows $f$ not to be odd.

\begin{theorem}\label{thm:properties}
  Assume that $N \geq 1$, $f \in C(\mathbb{R}, \mathbb{R})$ satisfies $(\hyperlink{f1}{f1}) - (\hyperlink{f3}{f3})$, and in addition $f$ is locally Lipschitz continuous when $N = 1$. Let $m \succeq_f m^*$, where $m^* \geq 0$ is the number given by Theorem \ref{thm:minimizers}. Then any minimizer $v \in S_m$ of \eqref{eq:Inf_m} satisfies the following properties:
    \begin{itemize}
      \item[$(i)$]  $v$ has constant sign,
      \item[$(ii)$] $v$ is radially symmetric up to a translation in $\mathbb{R}^N$,
      \item[$(iii)$] $v$ is monotone with respect to the radial variable.
    \end{itemize}
\end{theorem}

To prove Theorem \ref{thm:properties}, we make use of methods of ordinary differential equations when $N = 1$ and  this requires the locally Lipschitz continuity. In the higher dimensional case $N \geq 2$, the argument is based on \cite[Theorem 2]{M09} and \cite[Lemma 3.2]{BZ88}.

 In \cite[Theorem 0.1]{HT19}, by developing an original approach based on the introduction of a free functional and without assuming the nonlinearity to be odd, the existence of an energy ground state obtained in \cite{Sh14} was recovered under the assumption that the infimum $E_m$ can be achieved by a radial function. Our Theorem \ref{thm:properties} proves that it is indeed the case.

Our last theorem is the heart of the present paper, it answers positively for nonlinear scalar field equations a long standing issue.  To explain what is at stake, we need the following definition.

\begin{definition}\label{def:LAS}
  For given $\mu > 0$, a nontrivial solution $w \in H^1(\mathbb{R}^N)$ of the free problem
    \begin{equation*}\tag{$Q_\mu$}\label{eq:Q_mu}
      \left\{
             \begin{aligned}
               - \Delta u & = f(u) - \mu u \quad \text{in}~\mathbb{R}^N,\\
               u & \in H^1(\mathbb{R}^N),
             \end{aligned}
      \right.
    \end{equation*}
  is said to be a least action solution if it achieves the infimum  of the $C^1$ action functional
    \begin{equation*}
      J_\mu(u) := I(u) + \frac{1}{2}\mu\int_{\mathbb{R}^N}|u|^2 dx
    \end{equation*}
  among all the nontrivial solutions, namely
    \begin{equation*}
      J_\mu(w) = A_\mu := \inf\{ J_\mu(u) ~ | ~ u \in H^1(\mathbb{R}^N)\setminus\{0\}, J'_\mu(u) = 0\}.
    \end{equation*}
  For future reference, the value $A_\mu$ is called the least action of \eqref{eq:Q_mu}.
\end{definition}

Following \cite{DST21}, we note that any critical point $u$ of $J_{\mu}$ is also a critical point of $I$ restricted to $S_m$ (where $m$ is the mass of $u$) and conversely, any constrained critical point $u \in S_m$ of $I$ is also a critical point of $J_{\mu}$ (where $\mu$ is the Lagrange multiplier of $u$). One may thus wonder if the notions of energy ground states and of least action solutions coincide, namely if a least action solution is an energy ground state and vice versa. It is now known, see for example \cite{CKS21,JL21,LRN20}, that under the assumptions $(\hyperlink{f1}{f1}) - (\hyperlink{f3}{f3})$ there may exist least action solutions which are not energy ground states.  Remains, however, the possibility that any energy ground state is a least action solution. This issue can make precise by stating the following problem.  \bigskip

{\bf Problem:}  We know that a ground state energy minimizer $v \in S_m$ is a nontrivial solution to \eqref{eq:Q_mu}, where $\mu = \mu(v) > 0$ is the associated Lagrange multiplier. Is the minimizer $v \in S_m$ a least action solution to  \eqref{eq:Q_mu}?
\bigskip

For some odd $f \in C(\mathbb{R}, \mathbb{R})$ satisfying $(\hyperlink{f1}{f1}) - (\hyperlink{f3}{f3})$ it is known that there exists a unique positive solution to \eqref{eq:Q_mu} and that it is a least action solution, see for example \cite{CS20,KOPV17,LRN20}. Thus, in such situations, Theorem \ref{thm:properties} implies that any energy ground state is a least action solution. However, apart in some particular cases of this type where some uniqueness property was used, the Problem remained undecided until recently.  A positive answer was given in \cite{FJMM20} for a related minimization problem associated to a biharmonic equation with a power nonlinearity, see \cite[Proposition 3.9 and Theorem 1.3]{FJMM20}. Also, very recently in \cite{DST21}, the authors answered positively the Problem assuming essentially that the function $t \mapsto f(t)/t$ is nondecreasing on $(0, \infty)$ which allows to benefit from a Nehari manifold in the sense given in \cite{SW10}, see \cite[Theorem 1.3]{DST21}  for more details.  Note that the results of \cite{DST21} also hold when the analog of problem \eqref{eq:Q_mu} is set on an arbitrary domain $\Omega \subset \mathbb{R}^N$.   Finally, we mention \cite{I21} in which the Problem was positively answered for a nonlinearity which is a sum of two powers.

Our result in that direction covers all the previous particular cases, at least when the associated equations are set on all the space $\mathbb{R}^N$.

\begin{theorem}\label{thm:LeastAction}
  Assume that $N \geq 1$, $f \in C(\mathbb{R}, \mathbb{R})$ satisfies $(\hyperlink{f1}{f1}) - (\hyperlink{f3}{f3})$, and in addition $f$ is locally Lipschitz continuous when $N = 1$. Let $m \succeq_f m^*$ and denote by $\mu(v)$ the Lagrange multiplier corresponding to an arbitrary minimizer $v \in S_m$ of \eqref{eq:Inf_m}, where $m^* \geq 0$ is the number given by Theorem \ref{thm:minimizers}. Then the following statements hold.
    \begin{itemize}
      \item[$(i)$] Any minimizer $v \in S_m$ of \eqref{eq:Inf_m} is a least action solution of \eqref{eq:Q_mu} with $\mu = \mu(v) > 0$. In particular,
            \begin{equation*}
              A_\mu = E_m + \frac{1}{2} \mu m.
            \end{equation*}
      \item[$(ii)$] For given $\mu \in \{\mu(v) ~ | ~ v \in S_m ~ \text{is a minimizer of} ~ \eqref{eq:Inf_m}\}$, any least action solution $w \in H^1(\mathbb{R}^N)$ of \eqref{eq:Q_mu} is a minimizer of \eqref{eq:Inf_m}, namely
            \begin{equation*}
              \|w\|^2_{L^2(\mathbb{R}^N)} = m \qquad \text{and}  \qquad I(w) = E_m.
            \end{equation*}
    \end{itemize}
\end{theorem}

\begin{remark}\label{rmk:comments}
  \begin{itemize}
    \item[$(i)$] The conclusions of Theorem \ref{thm:LeastAction} (ii) were also observed in \cite{DST21, FJMM20, I21} in the corresponding frames.
    \item[$(ii)$] For alternative variational characterizations of the energy ground states, in related problems, we refer to \cite{CGT21,CT20, HT19}. Note that in \cite{CGT21,CT20}, a variational characterization of the associated Lagrange multiplier is proposed, see also \cite[Theorem 1.2]{DST21} in that direction.
  \end{itemize}
\end{remark}

Apart that it clarifies the relation between the two notions of energy ground states and of least action solutions, Theorem \ref{thm:LeastAction} has various practical consequences. First, it shows that searching for energy ground states  does not permit to find new solutions to the problem \eqref{eq:Q_mu}. It also allows to transmit to energy ground states the properties that are known to hold for the least action solutions. As an example of this we refer to \cite[Theorem 1.4]{LeWe21} where this strategy, relying in \cite{LeWe21}  on \cite[Proposition 3.9]{FJMM20} which corresponds to our Theorem \ref{thm:LeastAction}, is put at work. Knowing that least action solutions are nonradial, for a certain range of $\mu > 0$, leads to the conclusion that energy ground states are nonradial when their $L^2$ norm is sufficiently small.

At the core of the proof of Theorem \ref{thm:LeastAction} is a mountain pass type characterization of the nontrivial solutions of \eqref{eq:Q_mu} with $\mu = \mu(v) > 0$, see Lemma \ref{lemma:MP} below. Such a characterization is of independent interest. Indeed, a related characterization about least action solutions proved to be very fruitful to control the possible loss of compactness at infinity in non-autonomous problems, see \cite{BJT08,JT03,JT03-ANS} in that direction. The present version, which highlights the role of the $L^2$ mass, seems to have not been formulated before.

\begin{lemma}\label{lemma:MP}
  Assume that $N \geq 1$, $\mu > 0$ and $f \in C(\mathbb{R}, \mathbb{R})$ satisfies the conditions $(\hyperlink{f1}{f1})$ and
	\begin{itemize}
      \item[$(f2)'$] \hypertarget{f2'}{}  when $N \geq 3$, $\limsup_{|t| \to \infty} \frac{|f(t)|}{|t|^{\frac{N + 2}{N - 2}}} < \infty$,  \newline
                     when $N = 2$, for any $\alpha > 0$
                       \begin{equation*}
                         \lim_{|t| \to \infty} \frac{f(t)}{e^{\alpha t^2}} = 0,
                       \end{equation*}
    \end{itemize}
  and suppose in addition that $f$ is locally Lipschitz continuous when $N = 1$. Then for any nontrivial critical point $w \in H^1(\mathbb{R}^N)$ of $J_\mu$, any $\delta > 0$ and any $M > 0$, there exist a constant $T = T(w, \delta, M) > 0$ and a continuous path $\gamma : [0, T] \to H^1(\mathbb{R}^N)$ satisfying
    \begin{itemize}
      \item[$(i)$] $\gamma(0) = 0$, $J_\mu (\gamma(T)) < -1$, $\max_{t \in [0, T]} J_\mu(\gamma(t)) = J_\mu(w)$;
      \item[$(ii)$] $\gamma(\tau) = w$ for some $\tau \in (0, T)$, and
        \begin{equation*}
          J_\mu(\gamma(t)) < J_\mu(w)
        \end{equation*}
      for any $t \in [0, T]$ such that $\|\gamma(t) - w\|_{H^1(\mathbb{R}^N)} \geq \delta$;
      \item[$(iii)$] $m(t) := \|\gamma(t)\|^2_{L^2(\mathbb{R}^N)}$ is a strictly increasing continuous function with $m(T) > M$.
    \end{itemize}
\end{lemma}

As we shall see the proof of Lemma \ref{lemma:MP} is rather direct when $N \geq 3$, relying ultimately on the presence of a Pohozaev identity associated to \eqref{eq:Q_mu}. The proofs when $N=1$ and $N=2$ require separate more technically elaborated arguments.

The paper is organized as follows. In Section \ref{sect:ENE} we prove Theorem \ref{thm:minimizers}. Section \ref{sect:SSM} is devoted to the proof of Theorem \ref{thm:properties}.  Finally, in Section \ref{sect:LAC} we prove Lemma \ref{lemma:MP} and Theorem \ref{thm:LeastAction}.

\medskip
{\bf Notations.} Throughout this paper, for any $p \in [1, \infty)$, $L^p(\mathbb{R}^N)$ is the usual Lebesgue space endowed with the norm
  \begin{equation*}
    \|u\|_{L^p(\mathbb{R}^N)} := \left(\int_{\mathbb{R}^N} |u|^p dx\right)^{1/p},
  \end{equation*}
and $H^1(\mathbb{R}^N)$ the usual Sobolev space endowed with the norm
  \begin{equation*}
    \|u\|_{H^1(\mathbb{R}^N)} := \left(\|\nabla u\|^2_{L^2(\mathbb{R}^N)} + \|u\|^2_{L^2(\mathbb{R}^N)}\right)^{1/2}.
  \end{equation*}
Moreover, for given $u \in H^1(\mathbb{R}^N)$ and any $s \in \mathbb{R}$, we define the scaling function
  \begin{equation*}
    s \diamond u := e^{Ns/2}u(e^s \cdot),
  \end{equation*}
which remains in $H^1(\mathbb{R}^N)$ and preserves the $L^2$ norm when $s \in \mathbb{R}$ varies.


\section{Existence and nonexistence}\label{sect:ENE}
This section aims to prove Theorem \ref{thm:minimizers} and in particular we shall show the existence and nonexistence of minimizers of \eqref{eq:Inf_m} for suitable range of the mass $m > 0$. As a necessary preparation, we have the following lemma the proof of which is standard.

\begin{lemma}\label{lemma:I}
  Assume that $N \geq 1$ and $f \in C(\mathbb{R}, \mathbb{R})$ satisfies $(\hyperlink{f1}{f1})-(\hyperlink{f2}{f2})$. Then the following statements hold.
    \begin{itemize}
		\item[$(i)$] For any bounded sequence $\{u_n\}$ in $H^1(\mathbb{R}^N)$,
                     \begin{equation*}
                       \lim_{n \to \infty} \int_{\mathbb{R}^N} F(u_n) dx = 0
                     \end{equation*}
                   if $\lim_{n \to \infty} \|u_n\|_{L^{\infty}(\mathbb{R}^N)} = 0$, and
                     \begin{equation*}
                       \limsup_{n \to \infty} \int_{\mathbb{R}^N} F(u_n) dx \leq 0
                     \end{equation*}
                   if $\lim_{n \to \infty} \|u_n\|_{L^{2 + 4/N}(\mathbb{R}^N)} = 0$.
      \item[$(ii)$] There exists $C = C(f, N, m) > 0$  such that
                      \begin{equation*}
                        I(u) \geq \frac{1}{4}\|\nabla u\|^2_{L^2(\mathbb{R}^N)} - C(f, N, m)
                      \end{equation*}
                    for any $u \in H^1(\mathbb{R}^N)$ satisfying $\|u\|^2_{L^2(\mathbb{R}^N)} \leq m$. In particular, $I$ is coercive on $S_m$.
    \end{itemize}
  \end{lemma}

To proceed further, we recall the global infimum
  \begin{equation*}
    E_m := \inf_{u \in S_m} I(u)
  \end{equation*}
and make below a detailed study of its basic properties.
\begin{lemma}\label{lemma:E_m}
  Assume that $N \geq 1$ and $f \in C(\mathbb{R}, \mathbb{R})$ satisfies $(\hyperlink{f1}{f1})-(\hyperlink{f3}{f3})$. Then the following statements hold.
    \begin{itemize}
      \item[$(i)$] $- \infty < E_m \leq 0$ for all $m > 0$.
      \item[$(ii)$] There exists $m_0 > 0$ such that $E_m < 0$ for any $m > m_0$.
      \item[$(iii)$] $E_m < 0$ for all $m > 0$ if \eqref{eq:key_1} holds, and $E_m = 0$ for small $m > 0$ if \eqref{eq:key_2} holds.
      \item[$(iv)$] For any $m > m' > 0$ one has
                      \begin{equation}\label{eq:Em}
                        E_m \leq \frac{m}{m'} E_{m'}.
                      \end{equation}
                    If $E_{m'}$ is achieved then the inequality is strict.
      \item[$(v)$] The function $m \mapsto E_m$ is nonincreasing and continuous.
    \end{itemize}
\end{lemma}
\proof $(i)$ By Lemma \ref{lemma:I} $(ii)$, $I$ is bounded from below on $S_m$ and thus $E_m > - \infty$. For a fixed $u \in S_m \cap L^{\infty}(\mathbb{R}^N)$, we have $\|\nabla (s \diamond u)\|_{L^2(\mathbb{R}^N)} \to 0$ and $\|s \diamond u\|_{L^{\infty}(\mathbb{R}^N)} \to 0$ as $s \to - \infty$. In view of Lemma \ref{lemma:I} $(i)$, $E_m \leq \lim_{s \to - \infty} I(s \diamond u) = 0$.

$(ii)$ From $(\hyperlink{f3}{f3})$ and Step 1 of \cite[Proof of Theorem 2]{Be83-1}, there exists $u \in H^1(\mathbb{R}^N)$ such that $\int_{\mathbb{R}^N} F(u) dx > 0$. For any $m > 0$, set $u_m := u(m^{- 1/N} \cdot \|u\|^{2/N}_{L^2(\mathbb{R}^N)} \cdot x) \in S_m$. Since
  \begin{equation*}
    \begin{split}
      I(u_m)
        & = \frac{1}{2} \int_{\mathbb{R}^N} |\nabla u_m|^2 dx - \int_{\mathbb{R}^N} F(u_m) dx\\
        & = \frac{m^{\frac{N-2}{N}}}{2 \|u\|^{2(N-2)/N}_{L^2(\mathbb{R}^N)}} \int_{\mathbb{R}^N} |\nabla u|^2 dx - \frac{m}{\|u\|^2_{L^2(\mathbb{R}^N)}} \int_{\mathbb{R}^N} F(u) dx\\
        & =: A m^{\frac{N-2}{N}} - B m =: g(m),
    \end{split}
  \end{equation*}
it follows that $E_m \leq I (u_m) = g(m) < 0$ for any sufficiently large $m > 0$.

$(iii)$ When \eqref{eq:key_1} holds, we choose $u \in S_m \cap L^{\infty}(\mathbb{R}^N)$. For
  \begin{equation*}
    D := \int_{\mathbb{R}^N} |\nabla u|^2 dx \bigg/ \int_{\mathbb{R}^N} |u|^{2 + 4/N} dx > 0,
  \end{equation*}
by \eqref{eq:key_1}, there exists $\delta > 0$ such that $F(t) \geq D |t|^{2 + 4/N}$ for all $|t| \leq \delta$. Since $\|s \diamond u\|_{L^\infty(\mathbb{R}^N)} \leq \delta$ for some $s < 0$, it is clear that
  \begin{equation*}
    \begin{split}
      E_m \leq I(s \diamond u)
        & \leq \frac{1}{2} \int_{\mathbb{R}^N} |\nabla (s \diamond u)|^2 dx - D \int_{\mathbb{R}^N} |s \diamond u|^{2 + 4/N} dx \\
        & = \frac{1}{2} e^{2s} \int_{\mathbb{R}^N} |\nabla u|^2 dx - D e^{2s} \int_{\mathbb{R}^N} |u|^{2 + 4/N} dx \\
        & = - \frac{1}{2} e^{2s} \int_{\mathbb{R}^N} |\nabla u|^2 dx < 0.
    \end{split}
  \end{equation*}

When \eqref{eq:key_2} is satisfied, there exists $C_f > 0$ such that $F(t) \leq C_f |t|^{2 + 4/N}$ for any $t \in \mathbb{R}$. By the Gagliardo-Nirenberg inequality,
  \begin{equation*}
    \int_{\mathbb{R}^N} F(u) dx \leq C_f C_N m^{2/N} \|\nabla u\|^2_{L^2(\mathbb{R}^N)} \qquad \text{for all}~u \in S_m.
  \end{equation*}
For any $m > 0$ small enough such that $C_f C_N m^{2/N} \leq 1/4$, we have
  \begin{equation*}
    I(u) \geq \frac{1}{4} \|\nabla u\|^2_{L^2(\mathbb{R}^N)} > 0,
  \end{equation*}
and thus $E_m \geq 0$. From Item $(i)$, it follows that $E_m = 0$ for $m > 0$ small.

$(iv)$ Let $t := m/m' > 1$. For any $\varepsilon > 0$ there exists $u \in S_{m'}$ such that $I(u) \leq E_{m'} + \varepsilon$. Clearly, $v := u(t^{-1/N} \cdot ) \in S_m$ and then
  \begin{equation}\label{eq:Em_inequality}
    \begin{split}
      E_m \leq I(v)
          & = t I(u) + \frac{1}{2} t^{\frac{N-2}{N}} \left(1 - t^{\frac{2}{N}}\right) \int_{\mathbb{R}^N} |\nabla u|^2 dx \\
          & < t I(u) \\
          & \leq \frac{m}{m'}(E_{m'} + \varepsilon).
    \end{split}
  \end{equation}
Since $\varepsilon > 0$ is arbitrary, we see that the inequality \eqref{eq:Em} holds. If $E_{m'}$ is achieved, for example, at some $u \in S_{m'}$, then we can let $\varepsilon = 0$ in \eqref{eq:Em_inequality} and thus the strict inequality follows.

$(v)$ By Item $(i)$ and \eqref{eq:Em}, it is clear that $E_m$ is nonincreasing. To show the continuity, we define for given $u \in S_1$ and any $m > 0$ the real function
    \begin{equation*}
      \Phi_u(m) := \frac{1}{m} I(u(m^{- \frac{1}{N}} \cdot )) = \frac{1}{2}m^{- \frac{2}{N}} \|\nabla u\|^2_{L^2(\mathbb{R}^N)} -  \int_{\mathbb{R}^N} F(u) dx.
    \end{equation*}
  It is clear that
    \begin{equation*}
      \frac{E_m}{m} = \inf_{u \in S_1}\Phi_u(m).
    \end{equation*}
  Since $\Phi_u(m)$ is a concave function of $m^{- \frac{2}{N}}$, it follows that $E_m /m$ is continuous in $m > 0$ and so is $E_m$. \hfill $\square$

\bigskip
\noindent
{\bf Proof of Theorem \ref{thm:minimizers}.} We define
  \begin{equation*}
    m^* := \inf \{m > 0 ~ | ~ E_m < 0\}.
  \end{equation*}
It is easily seen from Lemma \ref{lemma:E_m} that $m^* \in [0, \infty)$,
  \begin{equation}\label{eq:Em_sharp}
    E_m = 0 \quad \text{if} ~ 0 < m \leq m^*, \qquad \qquad E_m < 0 \quad \text{when} ~ m > m^*;
  \end{equation}
in particular, $m^* = 0$ if \eqref{eq:key_1} holds, and $m^* > 0$ if \eqref{eq:key_2} holds. Let us first show that if $0 < m < m^*$ then $E_m = 0$ is not achieved. Indeed, assuming by contradiction that $E_m = 0$ is achieved for some $m \in (0, m^*)$, we infer from Lemma \ref{lemma:E_m} $(iv)$ that
  \begin{equation*}
    E_{m^*} < \frac{m^*}{m} E_m = 0
  \end{equation*}
and this leads a contradiction since $E_{m^*} = 0$ by \eqref{eq:Em_sharp}. The rest is to prove that the global infimum $E_m$ is achieved when $m > m^*$.

Fix $m > m^*$ and let $\{u_n\} \subset S_m$ be any minimizing sequence with respect to $E_m$. It is clear that $\{u_n\}$ is bounded in $H^1(\mathbb{R}^N)$ by Lemma \ref{lemma:I} $(ii)$ and then one may assume that up to a subsequence $\lim_{n \to \infty} \int_{\mathbb{R}^N} |\nabla u_n|^2 dx$ and $\lim_{n \to \infty} \int_{\mathbb{R}^N} F(u_n) dx$ exist. Since $E_m < 0$ by \eqref{eq:Em_sharp}, we deduce that $\{u_n\}$ is non-vanishing, namely
  \begin{equation}\label{eq:nonvanishing}
    \lim_{n \to \infty} \left(\sup_{y \in \mathbb{R}^N} \int_{B(y, 1)} |u_n|^2 dx\right) > 0.
  \end{equation}
Indeed, if \eqref{eq:nonvanishing} were not true, then $u_n \to 0$ in $L^{2 + 4/N}(\mathbb{R}^N)$ by Lions Lemma \cite[Lemma I.1]{Lions84-2} and thus
  \begin{equation*}
    \lim_{n \to \infty} \int_{\mathbb{R}^N} F(u_n) dx \leq 0
  \end{equation*}
via Lemma \ref{lemma:I} $(i)$; noting that $I(u_n) \geq - \int_{\mathbb{R}^N} F(u_n) dx$, we obtain a contradiction:
  \begin{equation*}
    0 > E_m = \lim_{n \to \infty} I(u_n) \geq - \lim_{n \to \infty} \int_{\mathbb{R}^N} F(u_n) dx \geq 0.
  \end{equation*}
Since $\{u_n\}$ is non-vanishing, there exists a sequence $\{y_n\} \subset \mathbb{R}^N$ and a nontrivial element $v \in H^1(\mathbb{R}^N)$ such that up to a subsequence $u_n(\cdot + y_n) \rightharpoonup v$ in $H^1(\mathbb{R}^N)$ and $u_n(\cdot + y_n) \to v$ almost everywhere on $\mathbb{R}^N$. Set $m' := \|v\|^2_{L^2(\mathbb{R}^N)} \in (0, m]$ and $w_n := u_n(\cdot + y_n) - v$. It is clear that
  \begin{equation}\label{eq:L2}
    \lim_{n \to \infty} \|w_n\|^2_{L^2(\mathbb{R}^N)} = m - m'
  \end{equation}
and, using  the splitting result \cite[Lemma 2.6]{JL20},
  \begin{equation}\label{eq:splitting}
    E_m = \lim_{n \to \infty} I(u_n) = \lim_{n \to \infty} I(v + w_n) = I(v) + \lim_{n \to \infty} I(w_n).
  \end{equation}
We shall prove below a claim and then conclude the whole proof.

\smallskip
{\bf Claim.} $\lim_{n \to \infty} \|w_n\|_{L^2(\mathbb{R}^N)} = 0$. In particular, $m' = m$ by \eqref{eq:L2}.

Let $t_n := \|w_n\|^2_{L^2(\mathbb{R}^N)}$ for each $n \in \mathbb{N}^+$. If $\lim_{n \to \infty} t_n > 0$, then \eqref{eq:L2} gives that $m' \in (0, m)$. In view of the definition of $E_{t_n}$ and Lemma \ref{lemma:E_m} $(v)$, we obtain
  \begin{equation*}
    \lim_{n \to \infty} I(w_n) \geq \lim_{n \to \infty} E_{t_n} = E_{m - m'}.
  \end{equation*}
From \eqref{eq:splitting} and \eqref{eq:Em}, it follows
  \begin{equation*}
    E_m \geq I(v)  + E_{m - m'} \geq E_{m'} + E_{m - m'} \geq \frac{m'}{m} E_m + \frac{m - m'}{m} E_m = E_m.
  \end{equation*}
Thus necessarily $I(v) = E_{m'}$  and this shows that $E_{m'}$ is achieved at $v \in S_{m'}$. But then still using \eqref{eq:splitting} and \eqref{eq:Em}, we obtain a contradiction:
  \begin{equation*}
    E_m \geq E_{m'} + E_{m - m'} > \frac{m'}{m} E_m + \frac{m - m'}{m} E_m = E_m,
  \end{equation*}
and so the claim is proved.

\smallskip
{\bf Conclusion.} Clearly, $v \in S_m$ by the above claim and thus $I(v) \geq E_m$. Since the claim and Lemma \ref{lemma:I} $(i)$ imply that
  \begin{equation}\label{eq:w_n}
    \lim_{n \to \infty} \int_{\mathbb{R}^N} F(w_n) dx \leq 0,
  \end{equation}
we also have $\lim_{n \to \infty} I(w_n) \geq 0$. Therefore, by \eqref{eq:splitting} we get $E_m \geq I(v)$ and hence $E_m < 0$ is achieved at $v \in S_m$.   \hfill $\square$

\begin{remark}\label{rmk:compactness}
   One can deduce further that $u_n(\cdot + y_n) \to v$ in $H^1(\mathbb{R}^N)$. Indeed, from \eqref{eq:splitting}, \eqref{eq:w_n} and the fact that $I(v) = E_m$, it follows
    \begin{equation*}
      \|\nabla w_n\|_{L^2(\mathbb{R}^N)} \to 0 \qquad \text{as} ~ n \to \infty.
    \end{equation*}
  Since $\lim_{n \to \infty} \|w_n\|_{L^2(\mathbb{R}^N)} = 0$, we obtain the strong convergence. Besides, having obtained minimizers for any $m > m^*$, we can conclude a posteriori from Lemma \ref{lemma:E_m} $(iv)$ that
    \begin{equation*}
      E_{m + m'} < E_m + E_{m'} \qquad \text{for any} ~ m, m' > 0 ~ \text{with} ~ m + m' > m^*.
    \end{equation*}
  However, we do not know if the above strict inequalities can be proved a priori under our general conditions.
\end{remark}


\section{Sign, symmetry and monotonicity}\label{sect:SSM}

This section is devoted to the proof of Theorem \ref{thm:properties}. Unless otherwise noted, for given $\mu > 0$ we use the notations
  \begin{equation*}
    g_\mu(t) := - \mu t + f(t) \qquad \text{and} \qquad G_\mu(t) := - \frac{1}{2} \mu t^2 + F(t).
  \end{equation*}
To deal with the case $N = 1$, a special treatment is required. To be more precise, we shall make use of the following classification result which is deduced by means of simple methods of ordinary differential equations.
\begin{lemma}\label{lemma:1D}
  Assume that $N = 1$, $f$ is a locally Lipschitz continuous function on $\mathbb{R}$ satisfying $(\hyperlink{f1}{f1})$, and $w \in H^1(\mathbb{R})$ is a nontrivial critical point of $J_\mu$ for some $\mu > 0$. Then $w$ has a sign. More precisely we have
    \begin{itemize}
      \item[$(a)$] If $w$ is negative somewhere,   then $\zeta_- := \sup\{ t < 0 ~|~ G_\mu(t) = 0\} \in (- \infty, 0)$,
                     \begin{equation*}
                       g_\mu(\zeta_-) < 0,
                     \end{equation*}
                   and after a suitable translation of the origin $w$ satisfies
                     \begin{itemize}
                       \item[$(a1)$] $w(x) = w(- x)$ for any $x \in \mathbb{R}$,
                       \item[$(a2)$] $w(x) < 0$ for any $x \in \mathbb{R}$,
                       \item[$(a3)$] $w(0) = \zeta_-$,
                       \item[$(a4)$] $w'(x) > 0$ for any $x > 0$.
                     \end{itemize}
      \item[$(b)$] If $w$ is positive somewhere,   then $\zeta_+ := \inf\{ t > 0 ~|~ G_\mu(t) = 0\} \in (0, \infty)$,
                     \begin{equation*}
                       g_\mu(\zeta_+) > 0,
                     \end{equation*}
                   and after a suitable translation of the origin $w$ satisfies
                     \begin{itemize}
                       \item[$(b1)$] $w(x) = w(- x)$ for any $x \in \mathbb{R}$,
                       \item[$(b2)$] $w(x) > 0$ for any $x \in \mathbb{R}$,
                       \item[$(b3)$] $w(0) = \zeta_+$,
                       \item[$(b4)$] $w'(x) < 0$ for any $x > 0$.
                     \end{itemize}
    \end{itemize}
  In particular, $w$ is a translation of the unique solution to the initial value problem $- u'' = g_\mu(u)$ with $u(0) = \zeta_{-}$ (or $u(0) = \zeta_{+}$) and $u'(0) = 0$.
\end{lemma}
\proof By regularity $w \in C^2(\mathbb{R}, \mathbb{R})$ and thus
  \begin{equation}\label{eq:1D-solution}
    - w'' = g_\mu(w) \qquad \text{in} ~ \mathbb{R}.
  \end{equation}
Since $|w(x)|$ and $|w'(x)|$ decay to zero exponentially as $|x| \to \infty$, we have
  \begin{equation}\label{eq:1D-zero}
    \frac{1}{2}|w'(x)|^2 + G_\mu(w(x)) = 0 \qquad \text{for} ~ x \in \mathbb{R}.
  \end{equation}
Without loss of generality, we only consider the case when $w$ is negative somewhere. By translating the point where $w$ achieves its negative minimum  to the origin, one may assume that $w'(0) = 0$. In view of \eqref{eq:1D-zero},
  \begin{equation*}
    G_\mu(w(0)) = 0 \qquad \text{and} \qquad  \zeta_- > - \infty.
  \end{equation*}
Since $(\hyperlink{f1}{f1})$ gives that $G_\mu(t) < 0$ for any $t < 0$ close enough to the origin, we also have $\zeta_- < 0$. Now assume by contradiction that $g_\mu(\zeta_-) \geq 0$. Since $w(0) \leq \zeta_-$, there exists $x^* \in \mathbb{R}$ such that $w(x^*) = \zeta_-$. Then
  \begin{equation*}
    w'(x^*) = 0 \qquad \text{and} \qquad w''(x^*) = - g_\mu(\zeta_-) \leq 0
  \end{equation*}
via \eqref{eq:1D-zero} and \eqref{eq:1D-solution} respectively. If $g_\mu(\zeta_-) > 0$, then since whenever $w(x) = \zeta_-$ one also has $w'(x) = 0$ and $w''(x) < 0$, $w$ can never go above $\zeta_- < 0$, which is impossible. On the other hand, if $g_\mu(\zeta_-) = 0$, then by uniqueness the conditions
  \begin{equation*}
    w(x^*) = \zeta_- \qquad \text{and} \qquad w'(x^*) = w''(x^*) = 0
  \end{equation*}
imply $w \equiv \zeta_-$, which is also impossible. With the desired conclusion $g_\mu(\zeta_-) < 0$ at hand, there exists $\varepsilon > 0$ such that $G_\mu(t) > 0$ for any $t \in (\zeta_- - \varepsilon, \zeta_-)$. If $w(0) < \zeta_-$, then $w(x_*) \in (\zeta_- - \varepsilon, \zeta_-)$ for some $x_* \in \mathbb{R}$ and so
  \begin{equation*}
    \frac{1}{2}|w'(x_*)|^2 + G_\mu(w(x_*)) > 0.
  \end{equation*}
This contradicts \eqref{eq:1D-zero}, and therefore $w(0) = \zeta_-$. Since $w$ is the global solution of \eqref{eq:1D-solution} with the initial conditions $w(0) = \zeta_-$ and $w'(0) = 0$, the rest follows from a standard adaptation of some arguments in \cite[Proof of Theorem 5]{Be83-1}.  \hfill $\square$

\begin{remark}\label{rmk:1D}
  Even though the nonlinearity $f$ in Lemma \ref{lemma:1D} is locally Lipschitz continuous, the nontrivial critical points of $J_\mu$ (if exist) are not necessarily unique up to a translation in $\mathbb{R}$ and up to a sign since we allow $f$ to be not odd.
\end{remark}

In the higher dimensional case $N \geq 2$, the radial symmetry of minimizers will be obtained as a direct consequence of a general symmetry result in \cite{M09}, and the proof of the monotonicity relies on Lemma \ref{lemma:*} below. We remark that the first part of Lemma \ref{lemma:*} is well known and the second part is a simple corollary of \cite[Lemma 3.2]{BZ88}.
\begin{lemma}\label{lemma:*}
  Let $v$ be a nonnegative measurable function defined on $\mathbb{R}^N$ such that for any $\alpha > 0$ the function $(v - \alpha)^+$ belongs to $H^1(\mathbb{R}^N)$ and has compact support, and denote by $v^*$ the Schwarz rearrangement of $v$. Then
    \begin{equation}\label{eq:*}
      \int_{\mathbb{R}^N}|\nabla v^*|^2 dx \leq \int_{\mathbb{R}^N} |\nabla v|^2 dx.
    \end{equation}
  Moreover, if the equality in \eqref{eq:*} holds then the level set
    \begin{equation*}
      \chi_\alpha := \{x \in \mathbb{R}^N ~ | ~ v(x) > \alpha\}
    \end{equation*}
  is equivalent to a ball for any $\alpha \in (0, \text{ess sup}(v))$.
\end{lemma}

\medskip
\noindent
{\bf Proof of Theorem \ref{thm:properties}}  By Lemma \ref{lemma:1D}, the case $N = 1$ is proved. We treat below the case $N \geq 2$. For given minimizer $v \in S_m$ of \eqref{eq:Inf_m}, we set
  \begin{equation*}
    v^+ := \max \{0, v\} \qquad \text{and} \qquad v^- := \min\{0, v\}.
  \end{equation*}
If $m^\pm := \|v^\pm\|^2_{L^2(\mathbb{R}^N)} \neq 0$, then Lemma \ref{lemma:E_m} $(iv)$ gives that
  \begin{equation*}
    E_m = I(v) = I(v^+) + I(v^-) \geq E_{m^+} + E_{m^-} \geq \frac{m^+}{m} E_m + \frac{m^-}{m} E_m = E_m,
  \end{equation*}
and thus $E_{m^\pm}$ is achieved at $v^\pm \in S_{m^\pm}$. Using Lemma \ref{lemma:E_m} $(iv)$ again, we obtain a contradiction:
  \begin{equation*}
    E_m \geq E_{m^+} + E_{m^-} > \frac{m^+}{m} E_m + \frac{m^-}{m} E_m = E_m.
  \end{equation*}
Hence $v$ has constant sign. Since any minimizer of \eqref{eq:Inf_m} is a solution of \eqref{eq:Q_mu} for some $\mu > 0$ and then by regularity must be of class $C^1$, we also deduce from \cite[Theorem 2]{M09} that $v$ is radially symmetric up to a translation in $\mathbb{R}^N$. To proceed further, without loss of generality, we may assume that $v \geq 0$ and $v(x) = \overline{v}(|x|)$ for some one variable function $\overline{v} : [0, \infty) \to [0, \infty)$. By the fact that $v(x) \to 0 $ as $|x| \to \infty$, it can be seen that $v$ is bounded and for any $\alpha > 0$ the function $(v - \alpha)^+$ belongs to $H^1(\mathbb{R}^N)$ and has compact support. Since the Schwarz rearrangement $v^*$ satisfies
  \begin{equation*}
    v^* \in S_m \qquad \text{and} \qquad \int_{\mathbb{R}^N} F(v^*) dx = \int_{\mathbb{R}^N} F(v) dx,
  \end{equation*}
it follows from Lemma \ref{lemma:*} that $E_m \leq I(v^*) \leq I(v) = E_m$ and thus
  \begin{equation*}
    \int_{\mathbb{R}^N}|\nabla v^*|^2 dx = \int_{\mathbb{R}^N} |\nabla v|^2 dx.
  \end{equation*}
By Lemma \ref{lemma:*} again, for any $\alpha \in (0, \max(v))$, the level set $\chi_\alpha$ is equivalent to a ball. We now assume by contradiction that $\overline{v}$ is not nonincreasing. Then
  \begin{equation*}
    \overline{v}(r_2) > \overline{v}(r_1) > 0 \qquad \text{for some} ~ r_2 > r_1 \geq 0.
  \end{equation*}
Since $\overline{v}(r) \to 0$ as $r \to \infty$, there exists $r_3 > r_2$ such that $\overline{v}(r_3) = \overline{v}(r_1)$. Denoting $a := \overline{v}(r_1)$ and $b := \overline{v}(r_2)$, one may see that for any $\alpha \in (a, b)$ the level set $\chi_\alpha$ is nonempty but not equivalent to a ball. This gives a contradiction and thus $v$ is nonincreasing with respect to the radial variable.  \hfill $\square$


\section{Least action characterization}\label{sect:LAC}

In this section we give the proof of Lemma \ref{lemma:MP} and then use it to prove Theorem \ref{thm:LeastAction} in a unified way for all the dimensions.

\bigskip
\noindent
{\bf Proof of Lemma  \ref{lemma:MP}}
  When $N \geq 3$ and for the given $w \in H^1(\mathbb{R}^N)$, we set
  \begin{equation*}
    \gamma(t) := \left\{
      \begin{aligned}
        & w\left(\frac{\cdot}{t}\right), & & \qquad \text{for} ~ t > 0,\\
        & ~~~ 0, & & \qquad \text{for} ~ t = 0.
      \end{aligned}
    \right.
  \end{equation*}
Note that $m(t) := \|\gamma(t)\|^2_{L^2(\mathbb{R}^N)} = t^N \|w\|^2_{L^2(\mathbb{R}^N)}$ and by the Pohozaev identity
  \begin{equation*}
    \begin{split}
      J_\mu (\gamma(t))
        & = \frac{1}{2} t^{N - 2} \int_{\mathbb{R}^N} |\nabla w|^2 dx - t^N \int_{\mathbb{R}^N} G_\mu(w) dx \\
        & = \frac{1}{2} \Big(t^{N - 2} - \frac{N - 2}{N} t^N \Big) \int_{\mathbb{R}^N} |\nabla w|^2 dx.
    \end{split}
  \end{equation*}
Clearly, the function $J_\mu (\gamma(t))$ has a unique maximum at $t = 1$ and $J_\mu (\gamma(t)) \to - \infty$ as $t \to \infty$. Thus, for  any  $M > 0$ we can choose a large constant $T = T(w, M) > 0$ such that the continuous path $\gamma : [0, T] \to H^1(\mathbb{R}^N)$ satisfies, for any $\delta > 0$, Items $(i) - (iii)$ of Lemma \ref{lemma:MP}.

In the case of $N = 1$, without loss of generality, we only consider the situation when the given $w \in H^1(\mathbb{R})$ is negative somewhere. Then the statement $(a)$ of Lemma \ref{lemma:1D} holds and we can define a negative continuous function $W : \mathbb{R} \to \mathbb{R}$ by
  \begin{equation*}
    W(x) = \left\{
      \begin{aligned}
        & ~~ w(x), & & \qquad \text{for} ~ x \geq 0,\\
        & \zeta_- - x^4, & & \qquad \text{for} ~ x \in [- \varepsilon, 0), \\
        & \zeta_- - \varepsilon^4, & & \qquad \text{for} ~ x < -\varepsilon.
      \end{aligned}
    \right.
  \end{equation*}
Here $\varepsilon > 0$ is a chosen small constant such that
  \begin{equation}\label{eq:varepsilon}
    \frac{1}{2}|W'(x)|^2 - G_\mu(W(x)) = 8 x^6 - G_\mu(\zeta_- - x^4) < 0  \qquad \text{for} ~ x \in [-\varepsilon, 0),
  \end{equation}
and it follows from $G_\mu(\zeta_-) = 0$ and $g_\mu(\zeta_-) < 0$. Setting
  \begin{equation*}
    \gamma(t) := \left\{
      \begin{aligned}
        & W(|\cdot| - \ln t), & & \qquad \text{for} ~ t > 0,\\
        & ~~~~~~~~ 0, & & \qquad \text{for} ~ t = 0,
      \end{aligned}
    \right.
  \end{equation*}
one may see that the path $\gamma : [0, \infty) \to H^1(\mathbb{R})$ is continuous,
  \begin{equation*}
    m(t) := \|\gamma(t)\|^2_{L^2(\mathbb{R})} = \left\{
      \begin{aligned}
        & \|w\|^2_{L^2(\mathbb{R})} - \int^{- \ln t}_{\ln t} |w(x)|^2 dx, & & \qquad \text{for} ~ t \in (0, 1), \\
        & \|w\|^2_{L^2(\mathbb{R})}, & & \qquad \text{for} ~ t = 1, \\
        & \|w\|^2_{L^2(\mathbb{R})} + 2 \int^0_{- \ln t} |W(x)|^2 dx, & & \qquad \text{for} ~ t > 1,
      \end{aligned}
    \right.
  \end{equation*}
and
  \begin{equation*}
    J_\mu(\gamma(t)) = \left\{
      \begin{aligned}
        & J_\mu(w) - \int^{- \ln t}_{\ln t} \Big(\frac{1}{2} |w'(x)|^2 - G_\mu(w(x))\Big) dx, & & \qquad \text{for} ~ t \in (0, 1), \\
        & J_\mu(w), & & \qquad \text{for} ~ t = 1, \\
        & J_\mu(w) + 2 \int^0_{- \ln t} \Big(\frac{1}{2} |W'(x)|^2 - G_\mu(W(x))\Big) dx, & & \qquad \text{for} ~ t > 1.
      \end{aligned}
    \right.
  \end{equation*}
By the fact that $G_\mu(w(x)) < 0$ for $x \neq 0$ and \eqref{eq:varepsilon}, we have
  \begin{equation*}
    J_\mu(\gamma(t)) < J_\mu(w) \qquad \text{for} ~ t \neq 1
  \end{equation*}
and
  \begin{equation*}
      J_\mu(\gamma(t)) < J_\mu(w) - 2 G_\mu(\zeta_- - \varepsilon^4) \cdot (\ln t - \varepsilon)  \to  - \infty \qquad \text{as} ~ t \to \infty.
  \end{equation*}
Noting also that $m(t)$ is strictly increasing and $m(t) \to \infty$ as $t \to \infty$, for any $\delta > 0$ and $M > 0$ there exists a large constant $T = T(w, M) > 0$ such that $\gamma : [0, T] \to H^1(\mathbb{R})$ is a desired continuous path of Lemma \ref{lemma:MP} when $N = 1$.

To cope with the remaining case $N = 2$, we adapt some arguments from \cite[Proposition 2]{BJT08}. For the given $w \in H^1(\mathbb{R}^2)$, we define $\Psi : [0, \infty) \times (0, \infty) \to \mathbb{R}$ by
  \begin{equation*}
    \Psi(\theta, s) := J_\mu \left(\theta w\left(\frac{\cdot}{s}\right)\right) = \frac{1}{2} \theta^2 \int_{\mathbb{R}^2} |\nabla w|^2 dx - s^2 \int_{\mathbb{R}^2} G_\mu(\theta w) dx.
  \end{equation*}
It can be easily seen that
  \begin{equation*}
    \begin{split}
      \Psi_\theta(\theta, s) &= \theta \int_{\mathbb{R}^2} |\nabla w|^2 dx - s^2 \int_{\mathbb{R}^2} g_\mu(\theta w)w dx, \\
      \Psi_s(\theta, s) &= - 2s \int_{\mathbb{R}^2} G_\mu(\theta w) dx, \\
      \frac{d}{d \theta} \int_{\mathbb{R}^2} G_\mu(\theta w) dx &= \int_{\mathbb{R}^2} g_\mu(\theta w)w dx.
    \end{split}
  \end{equation*}
Since the Nehari and Pohozaev identities give respectively
  \begin{equation*}
    \int_{\mathbb{R}^2} g_\mu(w)w dx = \int_{\mathbb{R}^2} |\nabla w|^2 dx \qquad \text{and} \qquad \int_{\mathbb{R}^2} G_\mu(w) dx = 0,
  \end{equation*}
there exist two positive constants $\theta_1 < 1 < \theta_2$ such that
  \begin{equation*}
    \frac{d}{d \theta} \int_{\mathbb{R}^2} G_\mu(\theta w) dx > 0 \qquad \text{for} ~ \theta \in [\theta_1, \theta_2],
  \end{equation*}
and
  \begin{equation}\label{eq:G_mu}
    \int_{\mathbb{R}^2} G_\mu(\theta w) dx \left\{
      \begin{aligned}
        & < 0, & & \qquad \text{for} ~ \theta \in [\theta_1, 1),\\
        & = 0, & & \qquad \text{for} ~ \theta = 1, \\
        & > 0, & & \qquad \text{for} ~ \theta \in (1, \theta_2].
      \end{aligned}
    \right.
  \end{equation}
As a direct consequence,
  \begin{equation}\label{eq:Psi_s}
    \Psi_s(\theta, s) \left\{
      \begin{aligned}
        &  > 0, & &  \qquad \text{for} ~ (\theta, s) \in [\theta_1, 1) \times (0, \infty),\\
        &  = 0, & & \qquad \text{for} ~ (\theta, s) \in \{1\} \times (0, \infty),\\
        &  < 0, & & \qquad \text{for} ~ (\theta, s) \in (1, \theta_2] \times (0, \infty).
      \end{aligned}
    \right.
  \end{equation}
On the other hand, noting that
  \begin{equation*}
    \Psi_\theta(1, s) = \int_{\mathbb{R}^2} |\nabla w|^2 dx - s^2 \int_{\mathbb{R}^2} g_\mu(w)w dx = (1 - s^2) \int_{\mathbb{R}^2} |\nabla w|^2 dx,
  \end{equation*}
for any $s \neq 1$ there exists $\vartheta_s \in (0, 1)$ such that
  \begin{equation}\label{eq:Psi_theta}
    \Psi_\theta(\theta, s) \left\{
      \begin{aligned}
        &  > 0,& &  \qquad \text{for} ~ (\theta, s) \in [1 - \vartheta_s, 1 + \vartheta_s] \times (0, 1),\\
        &  < 0,& &  \qquad \text{for} ~ (\theta, s) \in [1 - \vartheta_s, 1 + \vartheta_s] \times (1, \infty).
      \end{aligned}
    \right.
  \end{equation}
Also, with at hand the continuous function
  \begin{equation*}
    h(t) :=
      \left\{
        \begin{aligned}
          &  \frac{g_\mu(t)}{t},& &  \qquad \text{for} ~ t \neq 0,\\
          &  - \mu,& &  \qquad \text{for} ~ t = 0,
        \end{aligned}
      \right.
  \end{equation*}
one may find a small constant $s^* \in (0, 1)$ such that
  \begin{equation}\label{eq:s*}
    \Psi_\theta(\theta, s) = \theta \left(\int_{\mathbb{R}^2} |\nabla w|^2 dx - s^2 \int_{\mathbb{R}^2} h(\theta w)w^2 dx \right) > 0 \qquad \text{for} ~ (\theta, s) \in (0, 1] \times (0, s^*].
  \end{equation}
Now for any $\delta > 0$ we fix a small constant $\varepsilon = \varepsilon(\delta) > 0$ such that $1 - \varepsilon > s^*$ and
  \begin{equation*}
    \left\|w\left(\frac{\cdot}{s}\right) - w \right\|_{H^1(\mathbb{R}^2)} < \delta \qquad \text{for} ~ s \in [1 - \varepsilon, 1 + \varepsilon],
  \end{equation*}
and denote by $\eta(t) = (\theta(t), s(t)) : [0, \infty) \to \mathbb{R}^2$ the piecewise linear curve joining
  \begin{equation*}
    \begin{split}
      (0, s^*) \to (1 - \theta^*, s^*) \to (1 - \theta^*, 1 - \varepsilon)
        & \to (1, 1 - \varepsilon)  \\
        & \to (1, 1) \\
        & \to (1, 1 + \varepsilon) \to (1 + \theta^*, 1 + \varepsilon) \to (1 + \theta^*, \infty).
    \end{split}
  \end{equation*}
Here $\theta^* = \theta^*(w, \delta) \in (0, 1)$ is a chosen constant satisfying
  \begin{equation*}
    \theta^* \leq \min\{1 - \theta_1, \theta_2 -1, \vartheta_{1 - \varepsilon}, \vartheta_{1 + \varepsilon}\},
  \end{equation*}
and each segment is horizontal or vertical. Let $0 =: t_0 < t_1 < \cdots < t_6 < t_7 := \infty$ be such that for each $k = 0, 1, \cdots, 7$, the element $\eta(t_k) \in \mathbb{R}^2$ is an end point of a linear segment of the piecewise linear curve $\eta$. We define
  \begin{equation*}
    \gamma(t) := \theta(t) w\left(\frac{\cdot}{s(t)}\right), \qquad t \geq 0.
  \end{equation*}
Then the function $J_\mu(\gamma(t)) = \Psi(\eta(t))$ is strictly increasing on $(t_0, t_1)$, $(t_1, t_2)$ and $(t_2, t_3)$ by \eqref{eq:s*}, \eqref{eq:Psi_s} and \eqref{eq:Psi_theta} respectively. One may also see that $J_\mu(\gamma(t))$ is constant on $(t_3, t_4)$, $(t_4, t_5)$ by \eqref{eq:Psi_s}, and strictly decreasing on $(t_5, t_6)$ and $(t_6, t_7)$ via \eqref{eq:Psi_theta} and \eqref{eq:Psi_s} respectively. Moreover, using \eqref{eq:G_mu},
  \begin{equation}\label{eq:-infty}
    \begin{split}
      J_\mu(\gamma(t))
        & = \frac{1}{2} (1 + \theta^*)^2 \int_{\mathbb{R}^2} |\nabla w|^2 dx - s^2(t) \int_{\mathbb{R}^2} G_\mu((1 + \theta^*) w) dx \\
        & \to - \infty \qquad \text{as} ~ t \to \infty.
    \end{split}
  \end{equation}
Finally we observe that the mass function $m(t) := \|\gamma(t)\|^2_{L^2(\mathbb{R}^2)} = \theta^2(t)s^2(t)\|w\|^2_{L^2(\mathbb{R}^2)}$ is strictly increasing and
  \begin{equation}\label{eq:+infty}
    m(t) = (1 + \theta^*)^2 s^2(t) \|w\|^2_{L^2(\mathbb{R}^2)} \to \infty \qquad \text{as} ~ t \to \infty.
  \end{equation}
Since for any $M > 0$ we can deduce from \eqref{eq:-infty} and \eqref{eq:+infty} the existence of a large constant  $T = T(w, \delta, M) > 0$ such that
  \begin{equation*}
    J_\mu (\gamma(T)) < -1 \qquad \text{and} \qquad m(T) > M,
  \end{equation*}
the continuous path $\gamma : [0, T] \to H^1(\mathbb{R}^2)$ is a desired one and this completes the proof of Lemma \ref{lemma:MP}.  \hfill $\square$

\bigskip
\noindent
{\bf Proof of Theorem \ref{thm:LeastAction}}. To prove Item $(i)$, denoting by $w \in H^1(\mathbb{R}^N)$ an arbitrary nontrivial critical point of $J_\mu$, we only need to show that
  \begin{equation*}
    J_\mu(w) \geq J_\mu(v) = E_m + \frac{1}{2} \mu m.
  \end{equation*}
For a fixed $\delta > 0$ and $M := m > 0$, let $\gamma : [0, T] \to H^1(\mathbb{R}^N)$ be the continuous path given by Lemma \ref{lemma:MP}. In view of Lemma \ref{lemma:MP} $(i)$ and $(iii)$, there exists $t_0 \in (0, T)$ such that
  \begin{equation*}
    \|\gamma(t_0)\|^2_{L^2(\mathbb{R}^N)} = m
  \end{equation*}
and thus
  \begin{equation*}
    \begin{split}
      J_\mu(w) = \max_{t \in [0, T]} J_\mu(\gamma(t))
        & \geq J_\mu(\gamma(t_0)) \\
        & = I(\gamma(t_0)) + \frac{1}{2} \mu \int_{\mathbb{R}^N} |\gamma(t_0)|^2 dx \\
        & \geq E_m + \frac{1}{2} \mu m.
    \end{split}
  \end{equation*}

We now prove Item $(ii)$. In view of Item $(i)$, an arbitrary least action solution $w \in H^1(\mathbb{R}^N)$ of \eqref{eq:Q_mu} satisfies
  \begin{equation}\label{eq:A_mu}
    J_\mu(w) = A_\mu = E_m + \frac{1}{2} \mu m.
  \end{equation}
Assume by contradiction that $\|w\|^2_{L^2(\mathbb{R}^N)} \neq m$. Then, for
  \begin{equation*}
    \delta := \left| \sqrt{m} - \|w\|_{L^2(\mathbb{R}^N)} \right| > 0 \qquad \text{and} \qquad M := m > 0,
  \end{equation*}
we have the continuous path $\gamma : [0, T] \to H^1(\mathbb{R}^N)$ given by Lemma \ref{lemma:MP}. Noting that by Lemma \ref{lemma:MP} $(iii)$ there exists $t_0 \in (0, T)$ such that
  \begin{equation*}
    \|\gamma(t_0)\|^2_{L^2(\mathbb{R}^N)} = m \qquad \text{and} \qquad \|\gamma(t_0) - w \|_{L^2(\mathbb{R}^N)} \geq \delta,
  \end{equation*}
it follows from Lemma \ref{lemma:MP} $(ii)$ a contradiction:
  \begin{equation*}
    \begin{split}
      J_\mu(w)
        & > J_\mu(\gamma(t_0)) \\
        & = I(\gamma(t_0)) + \frac{1}{2} \mu \int_{\mathbb{R}^N} |\gamma(t_0)|^2 dx \\
        & \geq E_m + \frac{1}{2} \mu m.
    \end{split}
  \end{equation*}
Since we have proved $\|w\|^2_{L^2(\mathbb{R}^N)} = m$, it is easy to see further that $I(w) = E_m$ by \eqref{eq:A_mu}. \hfill $\square$

\begin{remark}
  \begin{itemize}
    \item[$(i)$] Among the many properties of the global infimum $E_m$, only the definition itself and the fact that $E_m \leq 0$ is achieved when $m \succeq_f m^*$, are used in our proof of Theorem \ref{thm:LeastAction}.
    \item[$(ii)$] The existence of least action solutions to \eqref{eq:Q_mu} is never used in the proof of Theorem \ref{thm:LeastAction} $(i)$. In fact, one does not even know a priori whether \eqref{eq:Q_mu} admits a least action solution, as the nonlinearity $f$ allows to be not odd.
  \end{itemize}
\end{remark}


\section*{Acknowledgements}
\addcontentsline{toc}{section}{Acknowledgements}

The authors thank Prof. Thierry Cazenave for useful remarks on a preliminary version of this work. Sheng-Sen Lu acknowledges the support of the China Postdoctoral Science Foundation (No. 2020M680174) and the National Natural Science Foundation of China (Nos. 11771324 and 11831009).



\vspace{-0.11cm}

{
\small
\begin{thebibliography}{100}
\addcontentsline{toc}{section}{References}














\bibitem{Be83-1} H. Berestycki, P.L. Lions,
{\it Nonlinear scalar field equations I: Existence of a ground state},
Arch. Rat. Mech. Anal. \textbf{82} (1983) 313--346.




\bibitem{BZ88} J.E. Brothers, W.P. Ziemer,
{\it Minimal rearrangements of Sobolev functions},
J. Reine Angew. Math. \textbf{384} (1988) 153–-179.




\bibitem{BJT08} J. Byeon, L. Jeanjean, K. Tanaka,
{\it Standing waves for nonlinear Schr\"{o}dinger equations with a general nonlinearity: one and two dimensional cases},
Commun. Partial Differ. Equ. \textbf{33} (2008) 1113--1136.


\bibitem{CKS21} R. Carles, C. Klein, C. Sparber,
{\it On ground state (in-)stability in multi-dimensional cubic-quintic Schr\"{o}dinger equations},
arXiv:2012.11637, December 2020.
	
	
\bibitem{CS20}  R. Carles, C. Sparber,
{\it Orbital stability vs. scattering in the cubic-quintic  Schr\"{o}dinger equation},
Rev. Math. Phys., \textbf{33} (2021), Article number: 2150004.


\bibitem{C03} T. Cazenave,
{\it Semilinear Schr\"{o}dinger equations},
Courant Lecture Notes in Mathematics. \textbf{10}, American Mathematical Society, Providence, RI, 2003.


\bibitem{CL82} T. Cazenave, P.L. Lions,
{\it Orbital stability of standing waves for some nonlinear Schr\"{o}dinger equations},
Commun. Math. Phys. \textbf{85} (1982) 549--561.




\bibitem{CGT21} S. Cingolani, M. Gallo, K. Tanaka,
{\it Normalized solutions for fractional nonlinear scalar field equations via Lagrangian formulation},
Nonlinearity \textbf{34} (2021), 4017--4056.


\bibitem{CT20} S. Cingolani, K. Tanaka,
{\it Ground state solutions for the nonlinear Choquard equation with prescribed mass}, to appear on Geometric Properties for Parabolic and Elliptic PDE's, INdAM Springer Series, Cortona 2019.




\bibitem{DST21} S. Dovetta, E. Serra, P. Tilli,
{\it Action versus energy ground states in nonlinear Schr\"{o}dinger equations}, Math. Ann. (2022). https://doi.org/10.1007/s00208-022-02382-z


\bibitem{FJMM20} A. J. Fernandez, L. Jeanjean, R. Mandel, M. Maris,
{\it Non-homogeneous Gagliardo-Nirenberg inequalities in $\mathbb{R}^N$ and application to a biharmonic non-linear  Schr\"{o}dinger equation},
J. Diff. Equ. \textbf{330} (2022) 1--65.


\bibitem{HaSt04} H. Hajaiej, C.A. Stuart,
{\it On the variational approach to the stability of standing waves for the nonlinear Schr\"{o}dinger  equation},
 Adv. Nonlinear Stud. \textbf{4} (2004) 469--501.


\bibitem{HT19} J. Hirata, K. Tanaka,
{\it Nonlinear scalar field equations with $L^2$ constraint: mountain pass and symmetric mountain pass approaches},
 Adv. Nonlinear Stud. \textbf{19} (2019) 263--290.


\bibitem{Ik14} N. Ikoma,
{\it Compactness of minimizing sequences in nonlinear Schr\"{o}dinger systems under multiconstraint conditions},
 Adv. Nonlinear Stud. \textbf{14} (2014) 115--136.








\bibitem{I21} Y. Ilyasov,
{\it On orbital stability of the physical ground states of the NLS equations},
	arXiv:2103.16353, August 2021.


\bibitem{JL19} L. Jeanjean, S.-S. Lu,
{\it Nonradial normalized solutions for nonlinear scalar field equations}, Nonlinearity \textbf{32} (2019), 4942--4966.


\bibitem{JL20} L. Jeanjean, S.-S. Lu,
{\it A mass supercritical problem revisited},
 Calc. Var. Partial Differ. Equ. \textbf{59} (2020), Article number: 174.


\bibitem{JL21} L. Jeanjean, S.-S. Lu,
{\it Normalized solutions with positive energies for a coercive problem and application to the cubic-quintic nonlinear  Schr\"{o}dinger equation},  Math. Models Methods Appl. Sci. \textbf{32} (2022), 1557--1588.









\bibitem{JT03} L. Jeanjean, K. Tanaka,
{\it A remark on least energy solutions in $\mathbb{R}^N$},
Proc. Amer. Math. Soc. \textbf{131} (2003), 2399--2408.


\bibitem{JT03-ANS} L. Jeanjean, K. Tanaka,
{\it A note on a mountain pass characterization of least energy solutions},
 Adv. Nonlinear Stud. \textbf{3} (2003), 461–-471.



\bibitem{KOPV17} R. Killip, T. Oh, O. Pocovnicu, M. Visan,
{\it Solitons and scattering for the cubic-quintic nonlinear Schr\"{o}dinger equation on $\mathbb{R}^3$ },
Arch. Ration. Mech. Anal. \textbf{225} (2017) 469--548.



\bibitem{LeWe21} E. Lenzmann, T. Weth,
{\it Symmetry breaking for ground states of biharmonic NLS via Fourier extension estimates},
arXiv:2110.10782, October 2021.

\bibitem{LRN20} M. Lewin, S. Rota Nodari,
{\it The double-power nonlinear Schr\"{o}dinger equation and its generalizations: uniqueness, non-degeneracy and applications},
Calc. Var. Partial Differ. Equ. \textbf{59} (2020), Article number: 197.






\bibitem{Lions84-1} P.-L. Lions,
{\it The concentration-compactness principle in the calculus of variations. The locally compact case, part 1},
Ann. Inst. H. Poincar\'{e} Anal. Non Lin\'{e}aire \textbf{1} (1984) 109--145.


\bibitem{Lions84-2} P.-L. Lions,
{\it The concentration-compactness principle in the calculus of variations. The locally compact case, part 2},
Ann. Inst. H. Poincar\'{e} Anal. Non Lin\'{e}aire \textbf{1} (1984) 223--283.






\bibitem{M09} M. Mari\c{s},
{\it On the symmetry of minimizers},
Arch. Ration. Mech. Anal. \textbf{192} (2009) 311--330.
















\bibitem{Sh14} M. Shibata,
{\it Stable standing waves of nonlinear Schr\"{o}dinger equations with a general nonlinear term},
Manuscripta Math. \textbf{143} (2014) 221--237.




\bibitem{St19} A. Stefanov,
{\it On the normalized ground states of second order PDE's with mixed power non-linearities},
Commun. Math. Phys. \textbf{369} (2019) 929--971.


\bibitem{St82} C.A. Stuart,
{\it Bifurcation for Dirichlet problems without eigenvalues},
Proc. Lond. Math. Soc. \textbf{45} (1982) 169--192.




\bibitem{SW10} A. Szulkin, T. Weth
{\it The method of Nehari manifold},  Handbook of nonconvex analysis and applications, Int. Press, Somerville, MA, (2010) 597--632.


\end{thebibliography}
}
\end{document}